\documentclass{amsart}
\usepackage{amsmath}
\usepackage{amssymb,latexsym}
\usepackage{graphicx}
\usepackage{epstopdf}
\usepackage{array}
\usepackage{tikz}
\usetikzlibrary{matrix}
\usepackage{color}
\newtheorem{theorem}{Theorem}
\newtheorem{definition}{Definition}

\newtheorem{proposition}{Proposition}
\newtheorem{lemma}{Lemma}
\newtheorem{corollary}{Corollary}
\newtheorem{remark}{Remark}
\def\T{\mathop{\mathcal{T}}}
\def\oT{\T\limits^\circ}

\def\opcup{\mathop{\cup}}
\def\opsqcup{\mathop{\sqcup}}

\begin{document}

 \title{On the decomposition of a 2D-complex germ with non-isolated singularities}

\author {N. C. Combe }

 \address{Aix-Marseille University, IML, Campus de Luminy case 907 13288 Marseille Cedex-9 France}
 
\email{noemie.combe@univ-amu.fr}

\thanks{ I am very grateful to A.Pichon for suggesting me the research subject, to N.A'Campo and D.B.Wajnryb for very stimulating discussions.}
\keywords{non-isolated singularities, Riemannian manifolds, Cheeger-Nagase metric, Hsiang-Pati metric }
\subjclass{53B20,14P10,30L05}
\date{}

\maketitle
\begin{abstract}
The decomposition of a two dimensional complex germ with non-isolated singularity, into semi-algebraic sets is given. This decomposition consists of four classes: Riemannian cones defined over a Seifert fibered manifold, a topological cone over thickened tori endowed with Cheeger-Nagase metric, a topological cone over
mapping torus endowed with Hsiang-Pati metric and a topological cone over the tubular neighbourhoods of the link's singularities.  In this decomposition there exist semi-algebraic sets that are metrically conical over the manifolds constituting the link. The germ is  reconstituted up to bi-Lipschitz equivalence to a model describing its geometric behavior. 
\end{abstract}

\section{Introduction}
The purpose of this paper is to study the structure of a complex 2-dimensional surface germ at the origin with one-dimensional singular locus  $(\Sigma,0)$. Our main result gives a splitting of $(X,0)$ in finite non-overlapping parts, except on the boundary. 
Let $(X,0)$ be a germ of surface singularities and $\mathcal{I}:(X,0)\to (\mathbb{C}^{3},0)$ an embedding. The embeding determines two  metrics. The so called  ``outer metric'',  which is the induced metric from the hermitian metric on $\mathbb{C}^{3}$ 
\[d_{\text{out}}(x_{1},x_{2})= | \mathcal{I}(x_{1})- \mathcal{I}(x_{2})|,\]
and the ``inner metric'' which is the arc length metric
\[
d_{\text{in}}(x_{1},x_{2}) = \inf_{\gamma \in \Gamma} \ell(\mathcal{I}\circ \gamma),\quad \Gamma =\{\gamma | \gamma: [0,1]\to X,\ \gamma(0)=x_{1}, \gamma(1)=x_{2}\},
\]
where $\ell(\mathcal{I}\circ \gamma)$ is the length of the arc in $\mathbb{C}^{3}$. 

 Let the surface germ $(X,0)$ be defined by a complex analytic germ $f:(\mathbb{C}^3,0) \to (\mathbb{C},0)$ and let  $\mathcal{B}_{\epsilon}$ be a ball,  for the induced metric, of small radius $\epsilon$ centered at the origin. Let $\mathcal{C}(X\cap S_{\epsilon})$ be  a real cone over the boundary of the ball intersected with the germ, where $S_{\epsilon}= \partial\mathcal{B}_{\epsilon}$ is the 5-sphere.
  
In the case of isolated singularities J. Milnor showed a conical structure lemma~\cite{Mi}. The lemma states that there exists a homeomorphism between the intersection $ X\cap \mathcal{B}_{\epsilon}$  and a real cone $\mathcal{C}(X\cap S_{\epsilon})$. Then the Milnor conical structure lemma states a topological equivalence between the two pairs $( \mathcal{B}_{\epsilon}, (X\cap \mathcal{B}_{\epsilon})$ and $(\mathcal{C}(S_{\epsilon}),\mathcal{C}(X)\cap S_{\epsilon}))$. 

In 1972 D. Burghelea and A. Verona  extended J. Milnor's conical structure lemma to the case of non-isolated singularities~\cite{BV} using the Whitney pre-stratification of a germ~\cite{Lo,Wh}. 
\begin{lemma}[ Conic structure lemma]\label{L: constlem}
{\sl Let $M$ be a complex three dimensional manifold, $X\subset M$ be a closed set and 
suppose that $\mathfrak{S}$ is a Whitney pre-stratification of $X$. Assume the strata of $\mathfrak{S}$ are connected. Consider a Riemannian metric on $M$ and let $d$ be the induced distance function.
 There exists a sufficiently small  $\epsilon_{0}$,  $0<\epsilon_{0}\ll1$
such that for any $\epsilon$, $0<\epsilon \leq\epsilon_{0}$, the pair $(\mathcal{B}_{\epsilon},X\cap\mathcal{B}_{\epsilon})$ is homeomorphic to the real cone on the pair $(S_{\epsilon},X\cap S_{\epsilon})$.} 
\end{lemma} 

The proof of this lemme used the Whitney pre-stratification of a germ~\cite{Lo,Wh}, and the fact that  any analytic set admits a Whitney pre-stratification. Let us recall the notion of the Whitney pre-stratification.

Let $M$ be a $n$-dimensional $C^{\mu}$-manifold without boundary and 
let $V_{1},V_{2} $ be a pair of  $C^{\mu}$-sub-manifolds of $M$, such that for each $p\in V_{i}$ there exists a coordinates chart $(\varphi_{i},U_{i})$ of class $C^{\mu}$ such that $p\in U_{i}$ and $\varphi_{i}(V_{i}\cap U_{i})=\mathbb{R}^{k}\cap \varphi_{i}(U_{i})$ for a suitable coordinate subspace $\mathbb{R}^{k}\subseteq \mathbb{R}^{n}$.

\vspace{5pt}
 \noindent {\bf Whitney condition~\cite{Ma,Wh}}. Let $r=dimV_{1}$ and let $ y\in V_{2}$. 
 The pair $(V_{1},V_{2}) $ satisfies the Whitney condition at  $y\in V_{2} $ if given any sequence $\{x_{i}\}$ of point of $V_{1}$ such that $x_{i}\to y $  and $\{y_{i}\}$ a sequence of points in $V_{2}$, also converging to $y$. Suppose  ${T V_{1}}_{x_{i}}$ converges to some r-plane $\tau \subset \mathbb{R}^n$ and that $x_{i} \neq y_{i}$  for all $i$ and the secant $\overline{(x_{i} y_{i})}$ to some $ \ell\subset \mathbb{R}^n$. Then $ \ell \subset \tau$.

The secant  $\overline{(x y)}$ refers to the line in $\mathbb{R}^{n}$ which is parallel to the line joining $x$ and $y\ne x$,  and passes through  the origin.

\begin{definition}[Whitney pre-stratification]
Let $M$ be a $C^{1}$ manifold and $ X\subset M$ be locally closed.
 A cover $\mathfrak{S}$ of $X$ by pairwise disjoint $C^{1}$ sub-manifolds is called a Whitney pre-stratification of $X$ if and only if :
 \begin{enumerate}
 \item It is locally finite. 
 \item If $V\in \mathfrak{S}$ then $(\overline{V}-V)\cap X$ is a union of strata. 
 \item $V_{i},V_{j} \in \mathfrak{S}$ satisfies the Whitney condition. 
\end{enumerate}
 \end{definition} 

Let $\Sigma$ be the singular locus and we define by induction $\Sigma^{j}=\Sigma(\Sigma^{j-1})$ be the j-th singular locus of the (j-1)-th-singular locus.
The complex algebraic surfaces $X$ defined by a polynomial $f$ have the following stratification : 
\[ \mathcal{V}_{1}:=X-\Sigma^{1},  \mathcal{V}_{2}:=\Sigma^{1}-\Sigma^{2},\dots , \mathcal{V}_{j}:=\Sigma^{j-1}-\Sigma^{j}, \dots ,  \mathcal{V}_{n}:=\Sigma^{n-1}-\Sigma^{n}.\]

A natural question in the study of the singularities, is to determine whether the homeomorphism defined in the conical structure lemma, $\varphi: \mathcal{M}= X\cap \mathcal{B}_{\epsilon}\to \mathcal{N}=\mathcal{C}(X\cap S_{\epsilon})$,  is   bi-Lipschitz, i.e. :
\[
\frac{1}{K} d_{ \mathcal{M}}(x_{1},x_{2}) \leq d_{ \mathcal{N}}(\varphi(x_{1}),\varphi(x_{2})\leq Kd_{ \mathcal{M}}(x_{1},x_{2}), \ K\geq 1.\]
The condition which assures  isomorphisms of metric spaces.
In most of the cases the answer to this question is negative. For weighted homogeneous isolated singularities only  $A_{1} $ and $D_{4}$ singularities have a bi-Lipschitz homeomorphism~\cite{BFN}. 

However, for isolated singularities it is possible to decompose a normal complex surface germ into a union of germs of pure dimension 4. These germs are partitioned into two distinct families~$({\large {\boldsymbol{\tau} },\tau })$~\cite{BM,BNP}, in a similar way of the decomposition  introduced by Margulis~\cite{KM} for symmetric space and extended by Thurston~\cite{Th} and Gromov~\cite{Gr2}.

 The so-called  thick part $\Large{\boldsymbol{\tau} }$ is such that the cone on its link $\mathcal{C}(\Large{\boldsymbol{\tau} } \cap S_{\epsilon})$ is bi-Lipschitz equivalent to $\Large{\boldsymbol{\tau} }\cap \mathcal{B}_{\epsilon}$ for the inner metric  on $X$.

The thin part { \Large$\tau $} of pure dimension 4 has a tangent cone of dimension strictly less than 4.  
Components of the thin part depend on the topology of their link. Non-trivial closed curves in the thin components have bounded lengths~\cite{Gr1, Pu}. As the curves  move to zero in the thin parts, their length reduces polynomially, while, for thick parts the length of closed curves reduces linearly.  
In the  case of  singular surfaces of the type $h(z_{3})=g(z_{1},z_{2})$, the classical tool to investigate the  isolated singularity is to use the carrousel~\cite{Le1,Le2,Le3} and the
Puiseux data~\cite{Le3,NS}. This Puiseux data has an incidence on the exponents defined in the metrics of the thin and thick parts.
We will briefly recall this procedure:
 let $\pi:(z_{1},z_{2},z_{3})\to (z_{1},z_{2})$ be a generic projection map and $\Gamma=\left\{\{h(z_{3})=g(z_{1},z_{2})\} \cap  \left\{\frac{\partial{h}}{\partial{z_{3}}}=0\right\}\right\}$ a curve that contains the singular locus of $X$.
 The discriminant locus is $\pi(\Gamma)=\Delta$ in the plane $\mathbb{C}^2$.
  The following definition and proposition are devoted to define the genericity of maps.
\begin{definition}
Let $(X,0)\subset (\mathbb{C}^n,0)$, with $dim(X)=p$.
Consider the plane in the set of Grassmannian $H\in\frak{G}(n-p,n)$ such that $0$ is an isolated point of the intersection between $H$ and $X$.
Then there exists an open set $U$of $0$ in $\mathbb{C}^n$ such that $X\cap H\cap U =\{0\}$ and such that the projection $\pi:X\cap U\to U'\subset H^{\bot}$ along $H$ is $k$-sheeted, $k\in\mathbb{N}$.
\end{definition}
\begin{proposition}
Let $f: (\mathbb{C}^3,0) \to \mathbb{C}$ be a complex germ and let $H\in\frak{G}(1,3)$ . 
The multiplicity of $\pi_{H}\vert_{_{ X}}$ is equal to $\text{ord}_{0}(f)$ if and only if $H\cap \mathcal{C}_{3}(X,0)=\{0\}$ where $\mathcal{C}_{3}(X,0)$ is the  Zariski tangent cone.
\end{proposition}
By a suitable change of coordinates we choose to have transversality of $z_{1}=0$ with the curves of $\Delta$ in $\mathbb{C}^2$.
Denote the equation of $\Delta$ by $q(z_{1},z_{2})=\prod_{i\geq1}q_{i}(z_{1},z_{2})$.

Let $q_{i}(z_{1},z_{2})$ be an irreducible component of $q(z_{1},z_{2})$. We consider the Puiseux developpement of $q_{i}^{-1}(0)$ given by:
\[\Phi_{p}^{i}(z)=\sum_{k\geq1} a_{k}^{i}z^{\eta_{k}}, \ a_{k}^{i}\in\mathbb{C},\ \eta_{k}\in \mathbb{Q}.\]

Let $v(q_{i})$ be the valuation of $q_{i}$. 
\begin{list}{-}{}
\item If  $q_{i}=0$ then $v(q_{i})=-\infty$.
\item If $\forall k, a_{k}^{i}\neq 0, k\in \mathbb{N}^{*}$ then $q_{i}$ does not have Puiseux pairs.
\item If the precedent cases do not hold, let $r_{1}=\inf\{ \eta_{ k} | a_{k}^{i}\neq 0 ,\eta_{ k}\in \mathbb{Q}\setminus\mathbb{N}\}$.
  The number $r_{1}=\frac{m_{1}}{n_{1}}$, with $gcd(m_{1},n_{1})=1$ defines the first Puiseux pair: $(m_{1},n_{1})$.
\item Let $r_{j+1}=\frac{m_{j+1}}{n_{1}\cdot...n_{j+1}}$ with $gcd(m_{j+1},n_{j+1})=1,\ n_{j+1}>1$. The pair $ (m_{j+1},n_{j+1})$ is the $(j+1)$-Puiseux pair of $q_{i}$.
\item The Puiseux characteristic exponents of $q_{i}$ are the rational numbers $r_{j}=\frac{m_{j}}{n_{1}\cdot...\cdot n_{j}}$
 where $1<r_{j}<r_{j+1}$ and $j\in\{1,2,...,s\}$. 
 \end{list}
 
The Puiseux characteristic data of $q$ is the union  of sets constituting the Puiseux characteristic data of each irreducible component $q_{i}$ of $q$.

\vspace{3pt}
The Puiseux exponents allows to describe a Waldhausen splitting of the Milnor fiber by a carousel diagram. Since we have Seifert-fibered spaces, we can have a decomposition as in Theorem 9.1.3 in~\cite{Wa}. Therefore, for a 2-dimensional compact oriented manifold $Y$, given a diffeomorphism $h$ on Y, one can define $h$ as the monodromy of a fibration over a circle. To define a 3-manifold fibered over a circle with monodromy $h$, we define the following mapping torus.
 
\begin{definition}[Mapping torus]
By mapping torus $MY$ over the manifold $Y$ we mean the quotient topological space :
\begin{equation}
MY:=([0,1]\times Y)/ \mathcal{R},
\end{equation}
 such that \begin{enumerate}
\item $ h:Y \rightarrow Y$ is an orientation preserving homeomorphism
\item $\mathcal{R}$
 is the equivalence relation given by : $ (1,x)\sim (0,h(x))$. 
\end{enumerate}
\end{definition}

In a case of germs with non-isolated singularities, the thick/thin decomposition   presents some difficulties. 
The difficulties are due to presence of singularities on the link $\mathcal{L}$ and, contrarily to the case of isolated singularities, the link $\mathcal{L}$ is not diffeomorphic to the boundary $\partial{F}$ of its Milnor fiber. 
In the case of non-isolated singularities we have follow two steps. 

\begin{enumerate} 
\item Firstly, it is necessary to make a normalization of the "non-normal" algebraic variety $X$. The normalization map  $n:\overline{X}\to X$ allows to introduce the  normal algebraic variety $\overline{X}$, which  has at most an isolated singularity. If $\overline{X}$ has an isolated singularity, we assume, without loss of generality, that the isolated singularity is at $0$.

 \item  Secondly, we perform a good resolution  $\rho:\widetilde{X}\to\overline{X}$ if $\overline{X}$  has an isolated singularity. By good resolution  we mean that exceptional divisor $\rho^{-1}(0)=\cup_{j} E_{j}$,
where $E_{j}\cong\mathbb{P}^{1}$ are the irreducible components and $ \rho^{-1}(0)$  has normal crossings. 
We assume that the proper transform of the singular locus $\Sigma$ of $X$ intersects $\cup_{i} E_{i}$ transversely, outside the intersections of the irreducible components $E_{i}$.

\end{enumerate}

Nevertheless, there exists some analogy between the boundary of the Milnor fiber $\partial{F}$ and the link $\mathcal{L}$.
Similarly to the decomposition proposed by Siersma~\cite{Si} of the boundary of  the Milnor fiber into two parts: $\partial{F}=\partial{F}_{1} \cup \partial{F}_{2}$, one can decompose the link  $\mathcal{L}=X\cap S_{\epsilon}$  into two sets $\mathcal{L}_{1}$,$\mathcal{L}_{2}$~\cite{NS},
\begin{equation}\mathcal{L}= \mathcal{L}_{1}\cup\mathcal{L}_{2}.\end{equation}

 The first set, $\mathcal{L}_{1}$ is the link without the interior of the tubular neighbourhoods of the singularities of the link. The second set, $\mathcal{L}_{2}$ is the intersection of these tubular neighbourhoods with the link $\mathcal{L}$. 
The set $\mathcal{L}_{1}$ is equivalent to the  boundary $\partial{F}_{1}$, in Siersma's decomposition of the Milnor fiber. The existence of a plumbed graph manifold for $\mathcal{L}_{1}$ is obtained by the Nemethi-Szilard algorithm~\cite{NS} for the boundray $\partial{F}_{1}$ and, hence, by  equivalence between $\partial{F}_{1}$ and $\mathcal{L}_{1}$, for $\mathcal{L}_{1}$.
To precise the set $\mathcal{L}_{2}$, we consider the decomposition of the singular locus of $X$ in irreducible components $\Sigma_{k}$:
\begin{equation}\label{E:singloc}
\Sigma=\bigcup_{k=1}^{n} \Sigma_{n}.
\end{equation}

Hence the singular locus of the link $\mathcal{L}$ is
  \begin{equation}{\Large\boldsymbol{\sigma}} =\bigcup_{k=1}^{n} \Sigma_{k}\cap \mathcal{L}.\end{equation}

 Let us  consider an irreducible component ${\Large\boldsymbol{\sigma}}_{k}=\Sigma_{k}\cap\mathcal{L}$ and  denote the tubular neighbourhood of ${\Large\boldsymbol{\sigma}}_{k}$ by $\mathcal{T}^{k}$. Observe that one can consider the tubular neighbourhood $\mathcal{T}^{k}$ as a total space with base space  ${\Large\boldsymbol{\sigma}}_{k}$ of a fibration. The fiber will be determined by
so-called transversal germs of ${\Sigma}_{k}$ defined at a point $p\in{\Sigma}_{k}$ with the condition that $p\neq0$.
The transversal germ 
at each non singular point $p\in{\Sigma}_{k}$, 
is obtained by intersecting $(X,0)$ by a linear space $\frak{H}$, transverse to ${\Sigma}_{k}$. The dimension of this linear space is : $dim(\frak{H})=codim(\Sigma_{k})$.
This transversal germ will be  denoted by:  $(\frak{H},p)$.
To define the fiber in the total space $\mathcal{T}^{k}$, one has to determine the Milnor fiber of $f$ restricted to the transversal slice $(\frak{H},p)$ , i.e: $\phi_{\epsilon, \delta}:B_{\epsilon}\cap f^{-1}|_{(\frak{H},p)}(\partial{B}_{\delta})\to\partial{B}_{\delta}$, for $\delta<\epsilon\ll 1$. This Milnor fiber does not depend on the point $p$ neither on the radius $\delta$. This Milnor fiber will be denoted by $M_{k}$.
Finally, the fibration: $\mathcal{T}^{k} \to \Large\boldsymbol{\sigma}_{k}$ has, up to homeomorphism equivalence, a cone over  $\partial{M}_{k}$ as fiber.

The aim is to study,  from a geometrical and topological point of view, the rate at which systoles, contained in the link $\mathcal{L}$ converge to $0$ in $X\cap \mathcal{B}_{\epsilon} $.
The answer to this question leads to a splitting of the neighbourhood of the origine: $X\cap \mathcal{B}_{\epsilon} $. The subspaces in $X\cap \mathcal{B}_{\epsilon} $ are Riemannian manifolds with metrics determining the rate of convergence of systoles.
These Riemannian manifolds have, topologically, the aspect of topological cones glued together by isometry.
The proof will rely on the conical structure lemma.

  \section{The decomposition model of $(X,0)$}

  In this section, the singular spaces $(X,0)$ will be described in terms of a model. This model provides a decomposition of  $X \cap \mathcal{B}_{\epsilon}$ and enables to study the rates at which systoles in the boundary $X \cap \mathcal{B}_{\epsilon}$ converge to zero.
  These Riemannian manifolds are topological cones endowed with Riemannian metrics. 

\begin{theorem}[Splitting thorem]\label{T:splitting}
Let $(X,0)\subset (\mathbb{C}^{3},0)$ be a two dimensional complex germ  defined by the form  $f:(\mathbb{C}^3,0)\to(\mathbb{C},0)$ with one-dimensional singular locus $ \Sigma $. Let us assume that  $0 \in \Sigma $.
Then there exists a splitting of  $X \cap \mathcal{B}_{\epsilon}$ into a finite family of compact and connected semi-algebraic sets $A_{i}$:
\begin{equation}
X \cap \mathcal{B}_{\epsilon}=\bigcup_{i=1}^{k}A_{i}.\end{equation}
Each of these semi-algebraic sets is one of the following type:
\begin{enumerate}
\item Riemannian cones defined over a Seifert fibered manifold with $0< t\leq \epsilon$.
\item Topological cones over thickened tori endowed with a Nagase metric.
\item Topological cones over mapping torus endowed with a Hsiang-Pati metric.
\item Topological cones over the tubular neighbourhoods of the link's singularities
\end{enumerate}
\end{theorem}
\begin{corollary}[Semi-algebraic partition]\label{C:semalgset}
Let $(X,0)\subset (\mathbb{C}^{3},0)$ be a two dimensional complex germ with a one-dimensional singular locus  $\Sigma$, with $0\in \Sigma$.
There exists a decomposition of  $(X,0)$ into a finite family of semi-algebraic sets, which can be:
\begin{enumerate}
\item Semi-algebraic germs which are a union of metrically conical sets, subsets in $X\cap \mathcal{B}_{\epsilon}$
\item Semi-algebraic germs of pure dimension 4 which have a tangent cone of dimension strictly less than $4$.
\end{enumerate}
\end{corollary}
\begin{remark}
This splitting  of which we show the existence, holds for any complex analytic surface germ with non-isolated singularities of complex dimension one. 
Somehow,  if the surface is of type $f(z_{1},z_{2},z_{3})= h(z_{3})-g(z_{1},z_{2})=0$, where $h$ is a polynomial function in $z_{3}$ and $g$ a polynomial function in variables $z_{1},z_{2}$, then the Puiseux exponents in the carousel decomposition play a role in the determination of the rates in the metrics.
\end{remark}
\begin{remark}
Let us come back to the notion of metrically conical. 
For semi-algebraic set  $\mathcal{J}$ in the germ $(X,0)$, 
consider a closed curve $\gamma$ in $\mathcal{J}$ which represents a non-trivial element of $\pi_{1}(\mathcal{J}-\{0\})$.  Define the action on the semi-algebraic set $\mathcal{J}$:   
\[\mathbb{R}_{+}^{*} \times \mathcal{J}\rightarrow \mathcal{J} :
(t,\gamma) \longmapsto t\gamma ,\] 
  restricted to $0\leq t\leq \epsilon$, with $\epsilon>0$ is small.

Let $h$ be a homeomorphism between  $\mathcal{J}\cap\mathcal{B}_{\epsilon}$ and the cone over the link $\mathcal{C}(\mathcal{J}\cap\mathcal{S}_{\epsilon})$. Suppose that the length of the curve $\gamma$ is of order $t^{\nu},\ \nu\geq 1$. If $h$ is a bi-Lipschitz map, then  $h(t\gamma)$ has also length of order $t^{\nu}$. 
On the other side the $\mathbb{R}_{+}^{*}$ action on the cone induces a bound, say $l$, on the length of the loop  $h(t\gamma)$.
Since the distance given by the bi-Lipschitz constant is $1/K$, therefore $ l\cdot h(t\gamma) $ is never closer to $0$ than $1/K$.
So, if this condition is not fulfilled, the map is not bi-Lipschitz and thus the semi-algebraic set is not metrically conical.
\end{remark}
 A topological cone  $\mathbb{R}_{+}^{*} \times \frak{M}$ over a Riemannian manifold $(\frak{M},g)$ equipped with the  metric : 
 \begin{equation}\label{E:metcon}
 g_{c}=dt^2+t^2g,
 \end{equation} is said  metrically conical.

The metrically conical set we  are interested in is a cone over the link $\frak{M} =X\cap\mathcal{S}_{\epsilon}$ with induced distance function bounded by $\epsilon>0$. 
A germ which is bi-Lipschitz equivalent to the cone over its link is called metrically conical.

\vspace{3pt}
Non-metrically conical algebraic sets have the following property:
  \begin{proposition}\label{P:DiamFib} Let $\Xi$ be the intersection of the non-metrically conical semi-algebraic set  with the sphere $\mathcal{S}_{\epsilon}$.
Then, the fiber of  $\Xi \to S^1$ has diameter of order $\alpha \epsilon^{\nu}$, where $\alpha>0$ is a constant, $\epsilon_{0}\geq \epsilon$ and $\nu>1$ is a rationnal number. 
 \end{proposition}
 
 \begin{corollary}\label{C:Sistol}
Let $A_{i}$ be a non-metrically conical semi-algebraic set, and let $0< t \leq \epsilon$.
If there exists a systole in the link of $A_{i}$, that is in $\Xi$, then this systole
 converges with order $O(t^{\nu})$ to $0$, as $t$ tends to 0.
 \end{corollary}
\begin{corollary}
Let $ A_{i}\subset\mathbb{R}\times\mathbb{R}^n$ and  $A_{i,t}:=A\cap(\{t\}\times\mathbb{R}^4)$.
For each non-metrically conical set $A_{i}$ there exists a Gromov-Hausdorff limit $\lim\limits_{t \to 0}(A_{t},d_{in}(A_{t}))$ where $d_{in}$ is the arc length metric.
\end{corollary}
Recall that a systole of a compact metric space $M$ is a metric invariant of $M$, defined to be the least length of a free loops representing nontrivial conjugacy classes in the fundamental group of $M$.

The proof follows from Bernig-Lytchak's Theorem 1 in~\cite{BL}.
Let $ X\subset\mathbb{R}\times\mathbb{R}^n$ be a 1-parameter family of subsets of $\mathbb{R}^4$. Suppose that each fiber $X_{t}:=X\cap(\{t\}\times\mathbb{R}^n)$ is connected. Then the Gromov-Hausdorff limit $\lim\limits_{t \to 0}(X_{t},d_{in}(X_{t}))$ exists.

\section{Proof of the decomposition theorem}\label{S:proof}

\subsection{Idea of proof}

To show that  $X\cap \mathcal{B}_{\epsilon}$ splits into  a finite family of semi-algebraic sets which can be classified along four types, we use the existence~\cite{BV} of a homeomorphism   $h:X\cap\mathcal{B}_{\epsilon}\to\mathcal{C}( X\cap \mathcal{S}_{\epsilon})$.  In the following part, we are interested in a small neighbourhood of $0$ in $X$. For the sake of simplicity, we denote by  $X$ the $\epsilon$-neighbourhood instead of $X\cap \mathcal{S}_{\epsilon})$.  Let us  give the draft of the proof:

\begin{enumerate} 
\item As was said before, we have to normalize $X$ and solve the singularity of the normalized surface $\overline{X}$ (if there exists an isolated singularity).
\item Construct a well chosen tubular neighbourhood of the exceptional divisor $\rho^{-1}(0)=\cup_{j} E_{j}$, where $E_{j}$ are the irreducible components and $ \rho^{-1}(0)$  has normal crossings.
 \item Apply Nagase and Hsiang-Pati procedure, on the tubular neighbourhood of the irreducible components of $\rho^{-1}(0)$:
 \begin{enumerate}
 \item On a neighbourhood of the irreducible components $E_{i}$, without the intersection points, following~\cite{HP}, there exists a Hsiang-Pati metric. For each neighbourhood of $E_{i}$, up to bi-Lipschitz, one can investigate the Hsiang-Pati metric on the  manifold $(0,1] \times [0,1]\times Y $, where $Y$ is a compact polygon in $\mathbb{R}^2$ with $\tilde{g}$ a standard metric: 
\begin{equation}\label{E:Hsiang-Pati } g_{_{HP}}= dt^2+t^2d\theta^{2}+t^{2\nu}\tilde{g}(y), \text{ where } \nu\geq1,\  (t,\theta,y)\in (0,1] \times [0,1]\times Y .
\end{equation}

\item Concerning the neighbourhood of the intersection points of the $E_{i}$, there exists a Cheeger-Nagase metric on $ (0,1]\times[0,1]^3$:
\begin{equation}\label{E:Nagase}
 g_{_{CN}}=dt^2+t^2d\theta^{2}+t^{2\nu}(ds^2+h^2(r,s)d\Theta^{2}),\ (t,\theta,s,\Theta)\in (0,1]\times[0,1]^3,
 \end{equation}
 where $h(r,s)$ is the product of smooth functions~\cite{Na}.
  \end{enumerate}
   
\item Finally, it will be necessary to go back to the link of $X$,  so as to characterize the splitting. The procedure allows to recover the topological cone over the link $\mathcal{L}$ of $X$  as a topological cone on the quotient of the normalization's link $\overline{\mathcal{L}}$ of $\overline{X}$ by some equivalence relation. 
    \end{enumerate}

\subsection{Local metrics in the components of $\overline{X}$ }

In this section, first we use a plumbing process to choose a tubular neighbourhood of the exceptional divisor in $\widetilde{X}$. An introduction of a continuous function, defined on the neighbourhood of $\rho^{-1}(0)$, allows to construct metrics.

Consider the normal surface $\overline{X}$ with isolated singularity at $0$.
Take into account that, for a surface with isolated singularity, the boundary of its Milnor fiber is diffeomorphic to its link $\overline{\mathcal{L}}$. 
A description of the boundary of the Milnor fiber can be obtained taking the boundary of the tubular neighbourhood of the exceptional divisor $\rho^{-1}(0)=\cup_{j} E_{j}\subset \widetilde{X}$ of $\overline{X}$. 
 More precisely, this tubular neighbourhood is the boundary of a 4-manifold, obtained by plumbing together the tubular neighbourhoods of each irreducible component $E_{j}$~\cite{BK,Ne,Or}. 
 
 Take the tubular neighbourhood $N_{i}$ of each irreducible components $E_{i}\cong\mathbb{P}^{1}$. This tubular neighbourhood is diffeomorphic to a disk-bundle  $(\frak{E}_{i},\frak{F}_{i},\frak p,D^{2})$ where $\frak{F}_{{i}}$ is the base space, $\frak{E}_{i}$ the oriented total space, $\frak p$ the continuous surjection satisfying local triviality conditions. In addition, we associate to this quadruplet its Euler number $e_{{i}}$.
 Note that $\frak{F}_{{i}}$ is a compact surface of genus $\frak g$, with boundary components.
For each pair $(E_{i},E_{j})$ of intersecting irreducible components, it is necessary to plumb the respective disk bundles. This procedure gives a tubular neighbourhood of $\rho^{-1}(0)$.
 Let us recall the plumbing process.
 Consider two irreducible components $E_{i},E_{j}$ that intersect at a point $q$ and consider their respective disk-bundles  
 $(\frak{E}_{i},\frak{F}_{i},\frak p,D^{2})$ and $(\frak{E}_{j},\frak{F}_{j},\frak p,D^{2})$. 
 
 Let $\mu:\frak{E_{i}}|_{D^{2}_{i}}\to\frak{E_{j}}|_{{D^{2}_{j}}}$, where the disks $D^{2}_{i}\subset \frak{F}_{i} $ and  $D^{2}_{j}\subset \frak{F}_{j}$, be a map between the restricted total spaces of these bundles.

Let us introduce two trivializations of the restricted bundles  $f_{i}:D^{2}_{i} \times D^{2}\to \frak{ E_{i}}|_{D_{i}^{2}},\,f_{j}:D^{2}_{j}\times D^{2}\to \frak{E_{j}}|_{D^{2}_{j}}$.
The way to obtain the plumbing of the two disk-bundles is to take the disjoint union of the total spaces $\frak{E}_{i}$ and $\frak{E}_{j}$ and to identify the images of $f_{i}$ and $f_{j}$.This identification is attained by pasting the component $D^{2}_{i} \times D^{2}$ of $\frak{E}_{i}$ to the component $D^{2}_{j} \times D^{2}$  of  $\frak{E}_{j}$ by the diffeomorphism $J=\left (\begin{smallmatrix} 0&1\\1&0\end{smallmatrix}\right):D^{2}_{i} \times D^{2}\to D^{2}_{j} \times D^{2}$. So, the map $\mu:\frak{E_{i}}|_{D^{2}_{i}}\to\frak{E_{j}}|_{{D^{2}_{j}}}$ is, in other words, defined  by the composition $\mu=f_{j}\circ J\circ f_{i}^{-1}:\frak{E_{i}}|_{D^{2}_{i}}\to\frak{E_{j}}|_{{D^{2}_{j}}}$.
By this procedure, the disk in the base space and the fiber are exchanged. So, this identification, defines the plumbing. Applying this procedure to all the intersection points of the irreducible components, one has the construction of a neighbourhood of the exceptional divisor  $\rho^{-1}(0)$.
Thus, referring to~\cite{Ne,Or} one can describe the link in terms of  graph manifold.  In the crudest sense, a graph manifold is a set of elementary pieces homeomorphic to $ S^{1} \times D^{2} $ or to $S^{1} \times (S^{2}\setminus\opsqcup_{i=1}^{n} D^2_{i}) $ that are glued together along some boundary components~\cite{Bo}.
Moreover, it is worth stressing that to each  graph manifold there exists a corresponding plumbing diagram of the plumbed  manifold. This plumbing diagram is a weighted graph. The vertices correspond to the Seifert-fibered manifolds, the edges correspond to the gluing tori and the weights are given by the Euler numbers $e_{{i}}$. 

The plumbed  manifold is a four dimensional manifold with boundary. 
For an isolated singularity, the boundary of the  plumbed manifold is the graph manifold and defines the topology of the link~\cite{Ne}.
We will choose a closed neighbourhood of $\rho^{-1}(0)\subset\widetilde{X}$ such that its boundary is the boundary of the  graph manifold.
On this neighbourhood, we will define a function $\frak{R}$ with range $[0,+\infty)$. This function will be the key to define the metrics 

\begin{definition}[Nagase~\cite{Na}]\label{D;Nafunct}
Let $\frak{R}$ be a function defined on a neighbourhood $\mathcal{W} \subset \widetilde{X}$ of $\rho^{-1}(0)=\opcup_{j}E_{j}$ and with range $[0,+\infty)$. We require that this function satisfies the following properties:
\begin{enumerate}{}{}
\item $\frak{R}|_{\rho^{-1}(0)}=0$
\item $\frak{R}|_{\mathcal{W}\setminus\rho^{-1}(0)}$ is smooth and positive
\item $\mathcal{W}\setminus\rho^{-1}(0)=\frak{R}^{-1}(0,1]$
\item $\frak{R}^{-1}(0,\epsilon]=(0,\epsilon]\times\frak{R}^{-1}(\epsilon)$
\end{enumerate}
\end{definition}
The last property is defined using $\frak{R}$ and appropriate flow lines in $\mathcal{W}\setminus\rho^{-1}(0)$ .
We split suitably  $\frak{R}^{-1}(\epsilon)=\opcup_{i} W_{i} \opcup_{j} V_{j}$ into  finite parts, which are non-overlapping except on the boundaries. 
These parts parts will be classified into two types: 
\begin{list}{-}{}
\item the tubular neighbourhood of an irreducible component $E_{i}$, outside the intersection points of the irreducible components, denoted by $W_{i}$. 
\item the neighbourhood, at the intersection points of the irreducible components, denoted by $V_{j}$.
 \end{list}
 Using Property 4 of Definition~\ref{D;Nafunct},  we split $\frak{R}^{-1}(0,\epsilon]=\opcup_{i} \frak{W}_{i}\cup_{j}\frak{V}_{j}$ into finite parts accordingly to the previous decomposition of $\frak{R}^{-1}(\epsilon) $ .
The part $ \frak{W}_{i}$ corresponds, in the splitting of $\frak{R}^{-1}(\epsilon)$, to the neighbourhood of the irreducible component  $E_{i}$ outside of the gluing tori. As for the part  $ \frak{V}_{j}$, it  corresponds to the gluing tori.
There is a bi-Lipschitz map from each $ \frak{W}_{j}$ with induced metric to a model $ (0,\epsilon]\times[0,1]\times Y$, where $Y$ is a compact polygon in $\mathbb{R}^2$. This polygon carries a Riemannian metric $\tilde{g}$. Moreover, this model has Hsiang-Pati metric ~\eqref{E:Hsiang-Pati }.
Whereas, each part $\frak{V}_{j}$ with induced metric is  bi-Lipschitz equivalent to $ (0,\epsilon]\times[0,1]^{3}$ with a Cheeger-Nagase metric~\eqref{E:Nagase}.
The exponents $\nu\geq 1$, in the metrics, correspond to the exponents in the divisor $\rho^{-1}(0)$.

 In the graph manifold, each Seifert-fibered manifold has a well defined Riemannian metric.
In the case where there exists a generic map $\pi:\mathbb{C}^{3}\to\mathbb{C}^{2}$ defined over $X$, one can construct a carousel decomposition. The Seifert-fibered manifolds in the link have topological data, given by the Puiseux exponents,  so it is the same for its metric. 
 Finally, for $\frak{W}_{i}$ the metric will depend on the Puiseux exponent, associated to $W_{j}$.
 For $\frak{V}_{i}$, notice that, according to L\^{e}'s description of the Waldhausen decomposition, by means of the  carousel~\cite{Le3}, the $\frak{V}_{i}$ need to be defined over thickened tori, i.e $\mathbb{T} \times I $ where $I=[a,b], b>a>0$ is an interval in $\mathbb{R}$ and $\mathbb{T}$ is a torus.
The real cone over the link $\mathcal{\overline{L}}$ results then, from gluing the boundary of the $\frak{W}_{j}$ with the boundary of the neighbouring  topological cones $\frak{V}_{i}$ and by using isometry so as to keep consistency of the metrics. The order in which they are glued  depends on the order in which the  manifolds $W_{j}$ and $V _{i}$ are glued in the plumbing diagram. %Moreover, it is necessary to use 

 A thickened torus is the Cartesian product of an annulus and a circle. 
 Let $A:=\{(r,\psi) | 1\leq r \leq 2,\  0\leq \psi \leq 2\pi  \} $, be an Euclidean annulus in polar coordinates, this annulus can be equipped of a family of metric $g_{a,b}= (b-a)^{2}dr^{2}+((r-1)b+(2-r)a)^{2}d\psi^{2}$ and for $o<a<b<1$ we can write $a=t^{\nu'}$and $b=t^{\nu}$ for some $t\in(0,1], \nu,\nu'>1$. We denote this metric by:
 \begin{equation}
g_{\nu,\nu'}^{t}:=(t^{\nu}-t^{\nu'})^2dr^2 +((r -1)t^{\nu}+(2-r)t^{\nu'})^2d\psi^2.
\end{equation} 
Notice that  $g_{\nu,\nu'}^{t}$ is isometric to the metric of an Euclidean annulus with inner radius $t^{\nu'}$ and outer radius $t^\nu$.
The metric $d t^2+t^2d\theta^2 +g_{\nu,\nu'}^{t}$ on $(0,1]\times {S}^1 \times A$ is bi-Lipschitz equivalent to a Cheeger-Nagase metric on $(0,1]\times[0,1]^3$,
\begin{equation}\label{E:HP}
\begin{aligned}
g_{_{HP}}&=d t^2+t^2d\theta^2 +g_{\nu,\nu'}^{t},\\
&=d t^2+t^2d\theta^2+ t^{2\nu}(ds^2+(s+t^{\nu-\nu'})^2d\psi^2), \quad s=(r-1)(1-t^{\nu'-\nu}), \\ 
&=d t^2+t^2d\theta^2+t^{2\nu}(ds^2+h(t,s)^2d\Theta^2), \quad h(t,s)= \frac{s+t^{\nu-\nu'}}{2\pi},\ \Theta \in [0,1],
\end{aligned}
\end{equation}
where, we call $g_{_{HP}}=d t^2+t^2d\theta^2 +g_{\nu,\nu'}^{t}$ the Cheeger-Nagase metric on $(0,1]\times {S}^1 \times A$.

Secondly, let us define the Hsiang-Pati metric on $\frak{W}_{j}$. 
The mapping torus $MY$ which was announced in the introduction, has a Riemannian metric ${g}_{\theta}$.

 Let $\theta \in [0,2\pi]$, then 
 for a small $\delta$:  \begin{equation}
{g_{\theta}}:= \begin{cases}  {g}_{0}, & \text{if } \theta \in [0,\delta] \\
 h^{*}{g}_{0} & \text{if }  \theta \in [2\pi-\delta, 2\pi] .\end{cases}
\end{equation}
Therefore, for any $t\in[0,1]$, the metric on $(0,1]\times MY$ is: 
  \begin{equation}\label{E:N}
  {g_{_N}}:=dt^2+t^2d\theta^2 +t^{2\nu}{g_{\theta}}, \quad t\in (0,1]. 
  \end{equation}
The compactification of the metric structure $(0,1]\times MY$ is done by  adding the singleton  $\{t=0\}$.

Finally,  we define a metric cone .
Let $(\frak{M},g)$ be a compact 3-manifold with a Riemannian metric $g$. A metric cone on $\frak{M}$ can be  obtained  by defintion of a metric on  $[0,1]\times \frak{M}$. This metric is the completion of the metric~\eqref{E:metcon} $ g_{_{c}}=dt^2+t^2g $ on $ t\in (0,1]\times \frak{M}$, by adding the singleton at $\{t=0\}$.  
\begin{remark}
Consider the mapping torus $MY$ over $ Y$.
Fix $t\in(0,1]$ and $\theta \in [0,2\pi]$. Then the diameter of $\{t\}\times \{\theta\} \times Y $ is of order $t^{\nu}$.
The rational number $\nu$ characterizes the rate at which $\{t\}\times \{\theta\} \times Y $ shrinks as $t$ tends to zero.
\end{remark}

\subsection{Local metrics on the subsets of $X$}
In this section we show that the metrics defined for $\mathcal{\bar{L}}$ hold for $\mathcal{L}$.
We assume that we have a plumbed graph obtained by the algorithm introduced by Nemethi and Szilard in~\cite{NS}. We will use the conditions showed in section 7.4 and Proposition 7.5.10, to show the splitting.

\begin{proposition}[Nemethi-Szilard]\label{T:NS}
Let  $\overline{\mathcal{L}}$ be the link of the normalized space.
Then, there exists an equivalence relation $\sim$ such that: $\mathcal{L}=\overline{\mathcal{L}}/\sim$.
\end{proposition}

A normalization of $X$ is a normal surface $\overline{X}$ and a holomorphic map $n:\overline{X} \to X$ such that :

\begin{enumerate} 
\item $n:\overline{X}\to X$ is proper and has finite fibers.
\item Let $\overline{\Sigma}=n^{-1}(\Sigma)$, then $\overline{X}-\overline{\Sigma}$ is dense in $\overline{X}$ and $n \arrowvert_{\overline{X}-\overline{\Sigma}}$ is a bi-holomorphic map. \end{enumerate}

Recall  that the normalized surface $\overline{X}$ is a disjoint union of irreducible components. So, the link $\mathcal{\overline{L}}$ of  $\overline{X}$ is a disjoint union of the links of each of these irreducible components.

Using the normalization map, for every $j$,  one can establish a correspondence between the strict transform of $\Sigma_{j}$
and a finite set $\{1,...,k\}$ of irreducible components in $\overline{X}$ which will be denoted $S_{j}^1,...,S_{j}^k$.

Notice that property (1) of the normalization implies that the degree of the map $n:S_{j}^i\to\Sigma_{j}$ is $d$.
On the other hand $\overline{S_{j,i}}:=S_{j}^{i} \cap\mathcal{\overline{L}}$ is diffeomorphic to a circle.
Therefore, for every $j\in\{1,2,...,n\}$, we define a degree $d$-map from $\overline{S_{j,i}}$ to a circle $S^{1}_{j}$:
\begin{equation}\label{E:10}
\mu_{i}:\overline{S_{j,i}}\to S^{1}_{j}\quad .
\end{equation}

Now we can define the equivalence relation as:
\begin{equation}\label{E:NS}
x\sim x' \iff \  \exists \ x\in \overline{S_{j,i}},\ x'\in\overline{S_{j,l}} \quad\text{ such that }\mu_{i}(x)=\mu_{l}(x'), \end{equation} 
For more ample explanations of 
Proposition~\ref{T:NS} we refer to Nemethi and Szilard~\cite{NS}.

Let us come back to the partition of $\overline{\mathcal{L}}$ into a union of manifolds $W_{i}$ and $V_{j}$,
$\overline{\mathcal{L}}=\opcup_{j} W_{j} \opcup_{i} V_{i}$. Recall that we have defined topological cones over each of the non-overlapping parts. These topological cones are labelled $\frak{W}_{i}$ and $\frak{V}_{j}$. The equivalence relation leaves invariant a subset of $\overline{\mathcal{L}}$. So, the topological cones  $\frak{W}_{i},\frak{V}_{j}$ that are defined over the manifolds of the link, invariant under the equivalence relation, are also   defined in  $\overline{X}$.

Recall, from introduction, that $\mathcal{L}=\mathcal{L}_{1}\cup\mathcal{L}_{2}$ where $\mathcal{L}_{1}:=\mathcal{L}\setminus\cup_{i=1}^{n}\oT{}^i$ and  $\oT{}^i$ is the interior of the tubular neighbourhood.
The second set $\mathcal{L}_{2}=(\cup_{i=1}^{n}\oT{}^i)\cap\mathcal{L}$ is the intersection of the tubular neighbourhoods $\oT{}^i$ with the link $\mathcal{L}$.
Thus we consider the cone $\mathcal{C}(\mathcal{L})=  \mathcal{C}(\mathcal{L}_{1}\cup\mathcal{L}_{2})$.
As was already said, $\T{}^i$ is a total space over $L_{i}=\Sigma_{i}\cap \mathcal{S}_{\epsilon}$ with fiber a cone over $\partial{M}_{i}$, where  $M_{i}$ is the Milnor fiber of $f$ restricted to $(\frak{h},p)$ for $\Sigma_{i}$.
Moreover, the tubular neighbourhood $\T{}^i$  has the homotopy type of $L_{i}$. 

Let us define the topological cone over the tubular neighbourhood of each singularity $L_{i}$ of the link $\mathcal{L}$. 

Note that from section 7.4 in Nemethi-Szilard~\cite{NS}, there exists a normal crossing divisor $\rho'^{-1}(\Sigma_{i})$.
Moreover, defining a small tubular neighbourhood $V$ of the discriminant curve for each point $p\in\Sigma_{i}, p\neq0$, the resolution $\rho'$ restricted to $\rho'^{-1}(V\cap(\{x\}\times \mathbb{C}^2))$ is an embedded resolution of the plane curve germ, at $p$.

Assume that the proper transform of $\Sigma_{i}$, intersects the union of irreducible components $\opcup_{j}{E}_{j}$ transversally and not in the intersection points of the irreducible components.
In order to have the link $\mathcal{L}$, it will be necessary to apply the equivalence relation in Nemethi-Szilard's proposition and glue together the circles $\overline{S_{i,j}}$ and $\overline{S_{i,l}}$, when the condition $\mu_{i}(\overline{S_{i,j}})=\mu_{l}( \overline{S_{i,l}})$ is fulfilled.
 The components $\overline{S_{i,j}}$ in the dual graph defined in section 7.4 of~\cite{NS}, 
are situated on the irreducible components $E_{r}$,  outside the intersection points of the irreducible components, in the resolution.
So,  $\overline{S_{i,j}}$  is on a Riemannian manifold with Hsiang-Pati metric. 
 Thus, the metric of the topological cone over the tubular neighbourhood $\T{}^i$ is of Hsiang-Pati type, with exponent $\nu$ depending on the exponent of the exceptional divisor.
Indeed, consider the resolution of the normalized space $\overline{X}$ and consider the irreducible component which contains the pre-images of the strict transform of 
$\Sigma_{i}$.
In order to define a type of Hsiang-Pati metric, take certain appropriate flow lines in the neighbourhood of $\overline{S_{i,k}}$. Proceed to the gluing of the circles $\overline{S_{i,k}}$ and $\overline{S_{i,l}}$ that have same image under $\mu $. Then, outside $S^{1}_{i}$ the type of Hsiang-Pati metric will be conserved on the topological cone defined over  the tubular neighbourhood $\T{}^i$in $\mathcal{L}$ .
The weight in the Hsiang-Pati metric of the topological cone defined over $\T{}^i$ is determined by the exponents in the resolution.

The boundary of this topological cone is glued to a boundary of the neighbouring Riemannian manifold: a topological cone over a thickened torus.
 Indeed, the neighbouring manifolds of the tubular neighbourhood $\T{}^i$ are gluing tori. 
Let $(\mathcal{C}_{top}(S^1 \times A),g_{_A})$ and $(\mathcal{C}_{top}( \T{}^i),g_{_{\T{}^i}})$ be Riemannian manifolds where $\mathcal{C}_{top}$ denotes a topological cone.
Let $f:\partial{C_{top}( \T{}^i)} \to \partial{\mathcal{C}_{top}({S}^1 \times A)}$ be a differentiable mapping and $\mathcal{C}_{top}( \T{}^i)$ be a Riemannian manifold with metric tensor $g_{_ {\T{}^i}}$.
The reciprocal image $f^{*}g_{_A}$ of $g_{_A}$ is a 2-covariant symmetric tensor on $\partial{\mathcal{C}_{top}(\T{}^i)}$ defined by: 
\begin{equation}
g_{_{\T{}^i}}=(f^{*}g_{_A})_{x}(v,w)=(g_{_A})_{y}(f'v,f'w)\  \forall v,w \in T_{x}\partial{\mathcal{C}_{top}(\T{}^i)}.\end{equation}

\subsection{The exponents in the metrics for $h(z_{3})=g(z_{1},z_{2})$ surfaces }
In this section we will use the well known method of carousel. This method will be applied to the case of the normalized surface $\overline{X}$. We give an outline of the method.

Let $\Delta$ be the discriminant curve in $\mathbb{C}^2$ and proceed to a Puiseux parametrization of each branch of $\Delta$. 
For each branch, it is necessary to truncate the Puiseux series in such a way that this truncation does not change the topology of the curve.
If one supposes that $\{z_{1}=\omega\}$, then one obtains a decomposition of the neighbourhood  of $\Delta$ in different regions, inversely proportional to the norm of the equation defining the singular locus.
For each tangent line to $\Delta$, $\{t= \alpha_{i} x\}$ we proceed as  follows. Suppose, without loss of generality, that there is only one tangent line $\{t= \alpha x\}$  to $\Delta$. 

Let us restrain the study of this curve to a set, centered at the tangent line: 
 \[A:=\{(x,y)| |x|\leq \epsilon , |y-\alpha x|\leq\mu|x|\}\subset \mathbb{C}^2.\]

Let $\nu_{r}$ be the biggest characteristic exponent in the truncated series.
For each $k\leq r$ let $\alpha_{k},\beta_{k},\gamma_{k}>0$ be positive constants, such that  $\alpha_{k} <\beta_{k}$.
 
The regions in $A$ are: 

\noindent -  $ \Upsilon_{k}:=\left\{(x,y):  \alpha_{k}|x^{\nu_{k}}\left|\leq |y-\sum\limits_{i=1}^{k-1}a_{i}x^{\nu_{i}}\right|\leq\beta_{k}|x^{\nu_{k}}|, \ \left|y-\sum\limits_{i=1}^{k-1}a_{i}x^{\nu_{k}}\right|\geq\gamma_{k}|x^{\nu_{k}}| \right\}, \ \forall k\in\{1,...,m\}$.\\

\noindent -  $\Omega_{i}:=\overline{ (\Upsilon_{i-1}-\Upsilon_{i})}, i\in\{ 2,..., m\}$.\\

\noindent -  $\Omega_{1}:=\overline{ (\partial A-\Upsilon_{1})}, i\in\{ 2,..., m\}$\\

\noindent - 
 $\Lambda:=\overline{ ( A-\cup \Omega_{i} \cup \Upsilon_{i}}), i\in\{ 1,..., q\}=\opcup\limits_{i=1}^{m+1}\Lambda_{i}$.\\

Consider the normalized surface $\overline{X}$ and $\pi:\overline{X}\to\mathbb{C}^2$ a generic projection with $(\pi(\Gamma))=\Delta \subset \mathbb{C}^2$.

The regions $\Omega,\Lambda,\Upsilon$ are lifted by $\pi$ to the boundary of the Milnor fiber of the surface.

Moreover, there exists a bi-Lipschitz map between the  regions $ \Omega_{i} , \Upsilon_{k},\Lambda$ and the topological cones defined on the gluing tori and the mapping tori.
\begin{itemize}
\item $\Omega_{i}$ is bi-Lipschitz equivalent to the topological cone over the gluing torus. The exponent in the metric of this topological cone is determined by the Puiseux exponents $\nu_{i}$ and $\nu_{i-1}$; 
\item $\Upsilon_{k}$ is bi-Lipschitz equivalent to the topological cone over the mapping torus defined over a planar surface $Y$ and with a finite order morphism $h$; The exponent in the metric is determined by the Puiseux exponent $\nu_{k}$.
\item $\Lambda$ is bi-Lipschitz equivalent to a topological cone over the mapping torus defined on the manifold $Y=D^2$. ~\cite{BNP}.
\end{itemize}

We will investigate the part of  $\overline{\mathcal{L}}$, which is invariant under the equivalence relation~\eqref{E:NS}. We recall that  this equivalence relation is defined such that:
 $\mathcal{L}=\overline{\mathcal{L}}/\sim$.
Define the topological cones, over the manifolds in $\overline{\mathcal{L}}$. 
Each topological cone carries a given metric with a given exponent $\nu$.
By carousel method, there exists a correlation between the zones in the carousel and the manifolds defining the link $\overline{\mathcal{L}}$. Moreover, a given Puiseux exponent gives a topological description of the corresponding manifold. 

Now take into consideration the topological cones with Puiseux exponents $\nu\geq1$, defined over the manifolds  of $\mathcal{\bar{L}}$ and invariant under the equivalence relation~\eqref{E:NS}.
Then, from the invariance of topological cones it follows that these manifolds are also defined on $\mathcal{L}$.
It is worth noticing that the invariant part of $\mathcal{L}$ is the set of plumbed manifolds equivalent to $\partial{F}_{1}$. In other words, the topological cones defined over the manifolds in $\mathcal{L}_{1}$ exist in the link $\mathcal{L}$. 

Take into account that if the Puiseux exponents of the pieces in $\mathcal{L}_{1}$  are greater than one, then the topological cones that are endowed with $\nu>1$ define non-metrically conical zones. In addition, if, in a manifold of $\mathcal{L}_{1}$, there exists a non-trivial loop, this loop will shrink at a rate greater than one.
 Thus we have defined a non-metrically conical zone in $X$.

As for the metrically conical zone, it corresponds to a Puiseux exponent $\nu=1$. 
So, the non-trivial loops defined in the Seifert manifold with rate 1, shrink linearly, that is $O(t)$ to 0.

First of all, one has the following non-commutative diagram between the discriminants of the normalized and non-normal surface:
\[
\begin{tikzpicture}
  \matrix (m) [matrix of math nodes,row sep=3em,column sep=4em,minimum width=2em]
  {
     \overline{X} & X \\
     \overline{\Delta} & \Delta \\};
  \path[-stealth]
    (m-1-1) edge node [left] {$\overline{\pi}$} (m-2-1)
            edge  node [above] {$n$} (m-1-2)
    (m-1-2) edge node [right] {$\pi$} (m-2-2)
          ;
\end{tikzpicture}
\]
Notice that the discriminant curves of $\overline{X}$ and $X$ differ by their Puiseux exponents. 
Indeed, not all of $\Delta$'s Puiseux data is preserved during the normalization process. 
In order to have a description of $\T{}^i$ we consider first $\partial{F}_{2,i}$.
Recall that for each $i$, $\partial{F}_{2,i}$ is a graph manifold. 
To show this assertion, it is sufficient to follow the Nemethi-Szilard algorithm~\cite{NS}. The decomposition into Seifert-fibered pieces of the three-manifolds $\partial{F}_{2,i}$ contains the data of $f$, restricted to the transversal slice germ of $\Sigma_{i}$.

\subsection{A property on the loops in $\partial{\mathcal{L}}_{2}^{i}$.}
Now, we will examine loops in the boundary $\partial{\mathcal{L}}_{2}^{i}$.

First, consider the set $\partial G_{\beta}:=\{ z \in \mathbb{C}^3 | |K(z)|^2=\beta \}$ where $K(z)$ generates the ideal defining the singular locus and let $\beta$ be a positive, small enough constant. 
Let us consider the map: 
\begin{equation}\label{E:fibre}\phi: (\mathcal{L}:=S_{\epsilon}\cap X)\cap \partial G_{\beta} \rightarrow S^{1}.\end{equation}
The fibers of this map are a set of $m$ manifolds. These manifolds will be denoted $R_{j}, j\in \{1,2,..., m\} $. 
Assume that the surface is such that the intersection of a hyperplane with the singular locus is transversal and contains 0.
Let $l$ be the projection operator onto the first coordinate axis $z_{1}$. Notice that in the case of  $l:\Sigma\to 0z_{1}$, one has a covering with finite multiplicity~\cite{Iom1}.
 Then, for a small $c\neq 0, l^{-1}(c)\cap\Sigma$ consists of a set of $m$ points that we will label $\{p_{1},p_{2},...,p_{m}\}$.
Consider the inverse image of these points $\{p_{i}\}_{i=1}^{m}$ in the surface. Then, we call $R_{j}$ the intersection of a sphere of small radius $\eta $, centered at a point $p_{j}$, and the pre-image of $c$ in the surface.

Going back to the previous map~\eqref{E:fibre}, if $\beta$ is sufficiently small then $S_{\epsilon}\cap G_{\beta}$ splits into a disjoint union of sets that are neighbourhoods of $S_{\epsilon}\cap\Sigma_{i}$.Thus one has to consider $\mathcal{L}^{i} :=(\mathcal{L}  \cap \partial{G_{\beta}^{i}})$ where $\partial{G_{\beta}}^{i}$ is a neighborhood of the $\Sigma_{i}$ component.
\begin{proposition}
If $i\in\{1,...,n\}$, then $\mathcal{L}^{i}$ stratifies over a circle with fiber $R_{i}$
\end{proposition}
An interesting topological information is given by I. Iomdin~\cite{Iom1}:
\begin{proposition}

$\forall z \in \mathcal{L}^{i}  \subset \mathbb{C}^3, i\in\{1,...,n\}$:
\[\lambda:\mathcal{L}^{i} \rightarrow V_{3,3}\ , \ z \mapsto [\text{grad f}, \text{grad g}, z]\]
 where
 $V_{3,3}$ is a complex Stiefel manifold.
%\end{proposition}

%\begin{proposition}
Let $ \lambda: \mathcal{L}^{i} \rightarrow V_{3,3}$ and let $i:R_{i} \rightarrow \mathcal{L}^{i} $ be an embedding.
Let $r_{j}$ be the fundamental cycle of the manifold $R_{i}$ and $s$ a generator of the first homology group $H_{1}(V_{3,3})$.
Then , \[\lambda_{*}\circ i_{*}(r_{i})=\kappa \cdot s, \] where $\kappa$ is an integer. 
\end{proposition}

So, using this procedure, we have defined a splitting of the germ $(X,0)$ into a finite family of semi-algebraic sets. Each of these sets correspond to a type of topological cone, endowed with Nagase, modified Nagase or Hsiang-Pati type of metric over the elements of the link of the singularity. 
Therfore this shows the existence of such a splitting of $X\cap\mathcal{B}_{\epsilon}$.
 \qed

\begin{proof} Propostion~\ref{P:DiamFib}.
Let $n: \overline{X} \to X$ be the normalization map and consider $X\setminus\Sigma$ with normalization $\overline{X}\setminus n^{-1}(\Sigma)$. Recall that outside the singular locus, the map is bi-holomorphic. So, outside the singular locus, the metric on $\overline{\mathcal{L}}\setminus n^{-1}(\Sigma\cap\mathcal{S}_{\epsilon})$ 
gives an equivalent metric, topologically speaking, on:  
$\mathcal{L}\setminus(\Sigma\cap\mathcal{S}_{\epsilon} )$. For normal surfaces, the boundary of its Milnor fiber is diffeomorphic to its link.
The discriminant curve of $\overline {X}$, provides the zones $\Omega_{k},\Upsilon_{i}, \Lambda$, delimited by the Puiseux-exponents of the discriminat curve. These zones give a filtration of the fibration:
$\mathcal{B}_{\epsilon}\cap f^{-1}(\partial{\mathcal{B}}_{\eta})\to \partial{\mathcal{B}}_{\eta},\ \eta<\epsilon,$ and provide a decomposition in terms of  graph manifolds of the link $\overline{\mathcal{L}}$.
We use the fact that there exists a bi-Lipschitz map between the set $\Upsilon_{k}$ with Puiseux exponent $\nu_{k}$, and the Reimannian manifold defined as a topological cone over a preimage $\pi^{-1}(\Upsilon_{k})$, with Riemannian metric $dt^2+t^2d\theta+t^{2\nu_{k}}g_{\theta}, t\in(0,1]$.  
 So, applying Cheeger-Kleiner theorem~\cite{CK}:  
let $\pi: \frak{W}_{j} \to \Upsilon_{j}$ be a Lipschitz map between metric spaces. We assume that for all $r\in\mathbb{R}^{+}$, if $U\subset \Upsilon_{j}$ and diam(U)$\leq r$, then the r-components of $\pi^{-1}(U)$ have diameter at most $Cr$.
It follows that if there exists  a systol in the manifold $MY$, which is a fiber of $\Xi \to S^{1}$ has diameter of order $\alpha t^{\nu}$, where $\alpha >0$ is a constant and $\nu$ is a Puiseux exponent. Note that $0<\alpha<\infty$ by section 3~\cite{BL}.
\end{proof}

\begin{proof}Corollary~\ref{C:Sistol}.
Let $\Xi$ be the region defined in the previous proposition. $\Xi$ is the region in $\mathcal{L}$ that corresponds to a zone of type $\Upsilon$, defined by a Puiseux exponent $\nu$. Thus, a systole, in this region, has a diameter of order $\alpha t^{\nu}$. Recall that the metric on the Riemannian manifold is of type : $dt^2+t^2d\theta+t^{2\nu_{k}}g_{\theta}, t\in(0,1]$. From a topological point of view, the topological cone can be studied in terms of a compact 1-parameter family of subsets of $ \mathbb{R}^n$. That is, if $A_{i}$ is the topological cone over $\Xi$, $A_{i}\subset\mathbb{R}\times\mathbb{R}^n$ and we define $A_{i,t}:=A_{i}\cap(\{t\}\times\mathbb{R})$. Notice that this set is connected.
So, as $t$ tends to $0$, the systole shrinks at a rate of order $O(t^{\nu})$.
\end{proof}

\subsection{Example}
Consider the complex germ defined by: 
\begin{equation}
  f(z_{1},z_{2},z_{3})=z_{3}^{2}-(z_{1}^3-z_{2}^2)^2(z_{1 }^4-z_{2}^3)=0,
   \end{equation}
 with singular locus $\Sigma=\{z_{3}=0, z_{1}^3-z_{2}^2=0\}$. 
The map $l:\Sigma \to z_{1}$ has finite multiplicity 2. So, the pre-image of $w\in\mathbb{C}^1$ is composed of two points $p_{1}$ and $p_{2}$.
Hence, the intersection in $X$: $l^{-1}(w)\cap X$, has an isolated singularity at $p_{1}$ and at $p_{2}$. 
Let us now define the following manifolds $R_{1}:=\mathcal{S}_{\epsilon}(c_{1})\cap l^{-1}(w)\cap X$ and $R_{2}:=\mathcal{S}_{\epsilon}(c_{2})\cap l^{-1}(w)\cap X$.
These two manifolds are diffeomorphic since, $p_{1}$ and $p_{2}$ are in the same branch of the singular curve defined by $\Sigma$.
We will denote this manifold for short $R$.
As it has been showed above, the set $X\cap\mathcal{S}_{\epsilon}\cap \partial{G_{\beta}}$ stratifies over a circle $S^1_{\gamma}$, $|z_{1}|=\gamma$ with fiber corresponding to the manifold $R$.

Let $n: \overline{X}\to X$ be the normalization map. The normalized space is defined by the following equation:
\begin{equation}\bar{f}=z_{3}^{2}-z_{1}^{4}-z_{2}^{3}=0.\end{equation}
 
By the Newton-Puiseux algorithm, one determines  the characteristic Puiseux exponent of the branch. In this example, the Puiseux characteristic exponent of the discriminant curve is 4/3. 
The delimited regions in this examples  are:
\begin{enumerate}
\item  let $\Upsilon_{1}:=\{( z_{1},z_{2}) \mathbb{C}^2:  \alpha |z_{1}^{\frac{4}{3} }|\leq |
z_{2}-z_{1}^\frac{4}{3}|\leq\beta|z_{1}^{\frac{4}{3}}| \}$, $\alpha\leq\beta\in \mathbb{R}$
\item $A:=\{(z_{1},z_{2})\in\mathbb{C}^2| |z_{1}|\leq \epsilon , |z_{2}^-z_{1}^|\leq\mu|z_{1}|\}\subset \mathbb{C}^2$
\item $\Omega_{1}:=\overline{ (\partial A-\Upsilon_{1})}, $
\item $\Lambda:=\overline{ ( A- \Omega_{1} \cup \Upsilon_{1}}), =\cup_{k=1}^{2}\Lambda_{k}$
\end{enumerate}

On the other side, applying Hirzebruch-Jung resolution of singularities, one obtains the resolution of this singularity and hence its dual representation graph. 
On this graph there exists a metrically conical zone, situated at the end vertex of the graph and a non-metrically zone, with metric rate $\nu=\frac{4}{3}$. 
In order to have the link of the non-isolated singularity, one must consider the inverse image of the singular locus $ \Sigma\subset X$ by $n$ in the link $\overline{\mathcal{L}}$. Then, apply the equivalence relation in Proposition~\ref{T:NS}.
In the neighbourhood of the glued points, we define a cone  over the boundary of $M_{j}$ and take the topological cone over this manifold. Recall that $M_{j}$ is the Milnor fiber of $f$ restricted to a transversal slice germ corresponding to $z_{1}^3-z_{2}^2$. This corresponds to the $\mathcal{L}_{2}$ part of the link. In $\mathcal{L}$, the adjacent manifolds are tori. These tori join $\mathcal{L}_{2}$ to the part $\mathcal{L}_{1}$. Notice that $\mathcal{L}_{1}$  is equivalent to the invariant part of $\overline{\mathcal{L}}$ under equivalence relation. Since the metric rate, for the non-metrically conical part of the link, is 4/3, we have the four types of topological cones in the splitting of $(X,0)$ :
 There exists a semi-algebraic set in $(X,0)$ which is metrically conical, a semi-algebraic set non-metrically conical with rate $4/3$ and semi-algebraic set equivalent to a topological cone  defined over $\mathcal{C}\partial{M_{j}}$.

\end{document}